\documentclass[11pt,reqno]{amsart}
\usepackage[a4paper,margin=30mm,headheight=14pt]{geometry}
\usepackage{amsmath,amssymb,amsthm}
\usepackage{xcolor,etoolbox,microtype}
\usepackage[T1]{fontenc}
\usepackage{lmodern}
\usepackage{hyperref}
\hypersetup{
  colorlinks=true,
  linkcolor=blue!45!black,
  citecolor=blue!45!black,
  urlcolor=blue!45!black,
  pdftitle={Smooth local isometric embeddings in dimension four},
  pdfauthor={Zhouyu Long},
  pdflang={en-US},
  pdfsubject={Smooth local isometric embeddings of four-dimensional Riemannian metrics},
  pdfkeywords={isometric embedding, Schlaefli problem, real principal type, Nash-Moser iteration, flat Gauss data}
}
\numberwithin{equation}{section}

\newtheorem{theorem}{Theorem}[section]
\newtheorem{proposition}[theorem]{Proposition}
\newcommand{\R}{\mathbb R}
\newcommand{\C}{\mathbb C}
\newcommand{\F}{\mathbb F}
\newcommand{\Sym}{\operatorname{Sym}}
\newcommand{\tr}{\operatorname{tr}}
\title[Local isometric embeddings in dimension four]{Smooth local isometric embeddings in dimension four}
\author[Zhouyu Long]{\bfseries Zhouyu Long$^{*}$}
\thanks{$^{*}$Department of Mathematical Sciences, Tsinghua University,
Beijing 100084, China.\newline
E-mail: \href{mailto:longzy25@mails.tsinghua.edu.cn}{\nolinkurl{longzy25@mails.tsinghua.edu.cn}}}

\makeatletter
\patchcmd{\@setauthors}{\centering\footnotesize}{\centering\large}{}%
  {\ClassError{amsart}{Unable to set author font}{Check the amsart author layout.}}
\patchcmd{\@setauthors}{\MakeUppercase{\authors}}{\authors}{}%
  {\ClassError{amsart}{Unable to preserve author case}{Check the amsart author layout.}}
\def\author@andify{\nxandlist{\unskip,\quad}{\unskip,\quad}{\unskip,\quad}}
\let\addresses\@empty
\makeatother
\date{}

\begin{document}
\begin{abstract}
Every smooth four-dimensional Riemannian metric admits a local
isometric embedding into $\R^{10}$, the optimal universal dimension.
A fixed-domain Nash--Moser argument prescribes any Gauss-compatible
rank-six second fundamental form with regular characteristic data.
For $n\ge4$, product-weighted cube evaluations give polynomial Gram
lifts of algebraic curvature tensors that are uniformly positive on a
fixed $n(n-1)/2$-dimensional support. In dimension four the support can
have a nonsingular complex characteristic quartic. With linear
tangential control, the transverse square-root estimate near affine
embeddings is sharp in $H^\ell$ for $\ell\ge1$.
\end{abstract}
\maketitle

\section{Introduction}

The theorem of Janet and Cartan \cite{Janet,Cartan} gives a local
isometric embedding of each real-analytic $n$-dimensional Riemannian
metric into $\R^{n(n+1)/2}$. We prove the smooth counterpart in
dimension four.

\begin{theorem}
\label{thm:main}
Let $(M,g)$ be a smooth four-dimensional Riemannian manifold and let
$p\in M$. There are a neighborhood $V$ of $p$ and a smooth embedding
$F:V\to\R^{10}$ such that $F^*\langle\cdot,\cdot\rangle=g$ on $V$.
The dimension $10$ is optimal for a universal smooth local theorem.
Moreover, in coordinates with $p=0$, $g(0)=I$, $dg(0)=0$, any six
independent normal forms $B$ satisfying the Gauss equation at zero
and having regular characteristic data can be prescribed, with
$F(0)=0$, $DF(0)=(I_4,0)^T$ and $D^2F(0)=(0,B)$.
\end{theorem}

Bryant, Griffiths and Yang proved characteristic regularity for general
second fundamental forms in dimension four
\cite[Theorem B and Corollary (1.c.6)]{BGY}; see also
\cite[Sections 3--4]{HanKhuri}. Goodman and Yang discuss
the smooth embedding application under a genericity condition on
the curvature jet $(R,\nabla R)$ \cite[Theorem 2.1(c)]{GoodmanYang}.
An earlier version \cite{LongV1} used an analytic flat reference
for the unrestricted smooth embedding theorem in dimension four.
Product-weighted cube evaluations give a common Gram support in every
$n\ge4$, and an explicit regular seed gives four-dimensional regularity.
Fixed polynomial kernel directions, the linear construction of
Nakamura--Maeda \cite{NakamuraMaeda}, and a finite-dimensional correction
give a single Nash--Moser application \cite[Section IV(b)]{BGY}
that also prescribes the second jet.
The domain and finite loss are fixed before rescaling the target.
Yang and Wang \cite[Theorems~1.1--1.2]{YangWang} study smooth local
isometric and conformal embeddings for $n\ge5$.

For any Euclidean-valued map $u$ of a coordinate domain, write
$\gamma(u)_{ij}=\langle \partial_i u,\partial_j u\rangle$.
For a four-dimensional immersion into $\R^{10}$, let
$B=(B_1,\ldots,B_6)$ denote the six normal forms in an orthonormal
normal frame, and let $B_{ij}$ be the vector with components
$(B_a)_{ij}$. We use the Gauss convention
\[
 \mathcal G(B)_{ijkl}
 =\sum_{a=1}^6\bigl((B_a)_{ik}(B_a)_{jl}
                    -(B_a)_{il}(B_a)_{jk}\bigr).
\]
The same convention applies to normal data in other dimensions,
with the sum taken over the normal directions.

When the six forms are independent, choose a basis $Q_1,\ldots,Q_4$
for their annihilator under the Frobenius pairing $A:B=\tr(A^TB)$.
With covectors represented in the given coordinates, set
\begin{equation}\label{eq:characteristic-quartic}
 M_B(\xi)_{i\cdot}=(Q_i\xi)^T,\qquad
 p_B(\xi)=\det M_B(\xi).
\end{equation}
We call the data regular if $p_B$ is nonzero as a polynomial and
$\partial_\xi p_B(\xi)\ne0$ at each nonzero real zero. A change
of annihilator basis multiplies $p_B$ by a nonzero scalar at each
base point. Thus the condition is independent of that basis and,
by the transformation of the principal symbol, of the coordinates.

We use $\nu\odot X=(\nu\otimes X+X\otimes\nu)/2$.

\section{Fixed-support Gram lifts}

\begin{theorem}\label{thm:gram-lifts}
For every $n\ge4$, put $E=\Sym_n$ and $r=n(n-1)/2$.
There is a fixed $r$-dimensional subspace $W\subset E$ and a
trace-one flat Gram tensor $S_*$ positive definite on $W$, such that
\[
 \mathcal G:\Sym^2 W\longrightarrow\mathcal R_n
 \quad\text{is onto},\qquad
 \dim\ker(\mathcal G|_{\Sym^2W})=\binom n4.
\]
Here $\mathcal R_n$ denotes the space of algebraic curvature tensors.
There is a polynomial choice $R\mapsto S_R$ in $\Sym^2 W$, with
$S_{R=0}=S_*$, such that on extending by zero on $W^\perp$,
\[
 S_R\ge\tfrac34 S_*,\quad \operatorname{rank}S_R=r,\quad
 \operatorname{im}S_R=W,\quad \mathcal G(S_R)=R .
\]
In dimension four $W$ can be chosen so that the characteristic
quartic $p_W$ has $\nabla p_W(\xi)\ne0$ for all nonzero $\xi\in\C^4$.
\end{theorem}

\begin{proof}
Identify $E$ with quadratic forms and its symmetric bilinear forms
with self-adjoint operators by the Frobenius pairing. The linear
Gauss map
$\mathcal G(T)_{ijkl}=T_{ik,jl}-T_{il,jk}$,
where $T_{ik,jl}=T(e_i\odot e_k,e_j\odot e_l)$, has kernel
$\Sym^4\R^n$ by the pair symmetries.
Put $m=n-1$, $v_s=(s,1)$ for $s\in\{-1,1\}^m$, and let $+$
denote the all-positive vertex. The degree-two evaluation space
has basis $1,s_i,s_is_j$ ($i<j$). Equip the cube with product
weights $\omega_s$ for which each coordinate is positive with
probability $p=1/(1+a)$, where $a=2m$. An orthonormal basis of
the evaluation space is $1,f_i,f_if_j$, with
$f_i=(1_{\{s_i=1\}}-p)/\sqrt{p(1-p)}$. Its reproducing kernel
$K$ satisfies
\[
 K(s,+)=\binom{k}{2}a^2+k(1-d)a+\frac{(d-1)(d-2)}2,
 \qquad k=\#\{s_i=1\},\quad d=m-k .
\]
Indeed each $f_i(s)f_i(+)$ is $a$ or $-1$. This expression is
positive for $k=0$. For $k=1$ it is $(d-1)((d-2)/2-a)<0$.
For $k\ge2$ it is positive since $a/2=m>d-1$ and the constant
term is nonnegative.
Subtract the rank-one evaluation term to form
\[
 S_\omega=\sum_s\omega_s v_s^{\otimes4}
              -K(+,+)^{-1}v_+^{\otimes4}.
\]
On the evaluation space this is the orthogonal projection
away from $K(\cdot,+)$. Hence it is positive semidefinite of rank $r$,
with kernel generated by the $m$ square differences $x_i^2-w^2$
and a quadratic $q_n$ evaluating to the nowhere-zero $K(\cdot,+)$.

Put $W_0=\operatorname{im}S_\omega$. The kernel quadrics
form a regular sequence: the square differences force any
common projective zero onto the cube, where $q_n$ never vanishes.
Their quotient thus has Hilbert series $(1+t)^n$ and dimension
$\binom n4$ in degree four. A fully symmetric tensor supported
on $W_0$ is a functional on its degree-four part, since its pairing on
quadrics is $(p,q)\mapsto\ell(pq)$. Therefore the supported Gauss
kernel has dimension $\binom n4=\dim\Sym^2 W_0-\dim\mathcal R_n$,
and the Gauss map is onto.

In dimension four this surjectivity persists for nearby planes.
A nearby support also carries a positive flat tensor, by transporting
$S_\omega$ and subtracting a continuous linear lift of its curvature.
For a nearby regular annihilator frame, interpolate the cube frame
with the regular frame of Appendix~\ref{app:regular-seed}.
The degree-nine gradient multiplication map is surjective at the
second endpoint by the proof there. Surjectivity is a Zariski-open
condition, so it holds for parameters arbitrarily close to the cube frame.

Choose $W$ as obtained, normalize the flat form to trace one,
and call it $S_*$. Pick a linear lift $T_R$ of $R$. On $W$ put
$D_R=S_*^{-1/2}T_RS_*^{-1/2}$ and
\begin{equation}\label{eq:gram-selection}
 S_R=T_R+(1+\operatorname{tr}(D_R^2))S_* .
\end{equation}
Congruence and $1+\lambda+\lambda^2\ge3/4$ give the asserted
uniform lower bound, and $\mathcal G(S_*)=0$ gives the curvature.
\end{proof}

In dimension four the positive Gauss fiber on this fixed $W$ is
exactly $\{T_R+sS_*:s>-\lambda_{\min}(D_R)\}$.
In every dimension the alternative analytic choice
$T_R+\sqrt{1+4\operatorname{tr}(D_R^2)}\,S_*$ has lower bound
$(\sqrt3/2)S_*$ and norm $O(1+|R|)$, by
$\sqrt{1+4s^2}-s\ge\sqrt3/2$ for $s\ge0$.

Blekherman's signed evaluation construction provides the quartic
precedent for the positive flat seed \cite[Theorem~2.7 and Section~6]{Blekherman}.
Yang and Wang \cite[Proposition~4.4]{YangWang} construct full-rank flat
data on a support spanned by rank-one quadrics. A basis of $r$ such
quadrics gives $r$ independent flat tensors in $\Sym^2 W$, by taking
the squares of the basis elements. In dimensions four and five this
exceeds $\binom n4$, the kernel dimension for a surjective
restricted Gauss map.

For $n\ge5$, the determinant of every complex $n\times n$
linear pencil in $n$ variables has a projective critical zero.
An injective pencil has a projective $(n-1)$-plane that meets the corank-at-least-two
determinantal variety of codimension four. A noninjective pencil
takes the zero matrix as its value at some nonzero point.
In both cases the determinant and its first derivatives vanish there.
For $n=2,3$ the dimension difference is zero, so a nonzero
flat tensor and surjectivity on a support of the stated rank
cannot coexist.

Appendices~\ref{app:regular-seed} and \ref{app:block-seed} give
explicit four-dimensional flat seeds. Appendix~\ref{app:analytic-jet}
realizes prescribed Gram data for analytic metrics by Cartan--K\"ahler.

\section{A fixed-domain tame right inverse}

\begin{proposition}\label{prop:metric-inverse}
Let $u_0$ be a smooth immersion of a neighborhood of $0\in\R^4$
into $\R^{10}$, with independent normal forms and regular
characteristic data at $0$, extended smoothly with compact support.
There are a fixed ball $U$ about $0$, a fixed spatial cutoff
$\vartheta=1$ near $\overline U$, and integers
$A_0\ge D$ such that each smooth background $u$ in a fixed
$H^{A_0}(\R^4)$ neighborhood of $u_0$ has a linear right inverse
$R_u$ on $U$ taking smooth Sobolev data $H$ to functions satisfying
\begin{equation}\label{eq:metric-tame}
 D\gamma_uR_uH=H\ \text{on }U,\qquad
 \|R_uH\|_{S-D}\le C_S(\|H\|_S+\|\vartheta u\|_S\|H\|_D),
 \quad S\ge D,\ S\in\mathbb Z .
\end{equation}
Here the norms are on $\R^4$ and one choice of inverse works
at all orders. For $S\ge A_0$, the map $(u,H)\mapsto R_uH$ extends
jointly continuously to pairs in $H^S\times H^S$ in the stated
background neighborhood, with values in $H^{S-D}$.
\end{proposition}

\begin{proof}
Localize $u$ by $\vartheta$, equal to one near the reference ball.
Choose orthogonal coordinates in which $dx_4$ is noncharacteristic at zero
and let $I=\{(i,j):1\le i\le j\le3\}$. Put $N_uV=(u_{ij}\cdot V)_{ij}$.
The matrix $M_u:V\mapsto((u_i\cdot V)_i,(N_uV)_I)$ is invertible
on a fixed small ball. Indeed, at the base point the annihilator
map $Q\mapsto Qe_4$ is invertible, so no nonzero annihilator form
is supported on the tangential slice. The six tangential normal
pairings are therefore independent. Write
$(C_0,C_1)=(N_u)_{I^c}M_u^{-1}$ and $h=H/2$. The identity
$\tfrac12D\gamma_uV=\operatorname{sym}\partial X-N_uV$, where
$X_i=u_i\cdot V$, gives the equivalent system
\begin{align*}
 L_cX&=(\operatorname{sym}\partial X)_{I^c}
         -C_1(\operatorname{sym}\partial X)_I-C_0X=h_{I^c}-C_1h_I=:f,\\
 V&=M_u^{-1}(X,(\operatorname{sym}\partial X)_I-h_I).
\end{align*}
Extend $c=(C_0,C_1)$ and $M_u^{-1}$ by fixed smooth cutoffs within
the frame-valid ball, equal to one on a smaller outer ball. The quotient $\pi A=A_{I^c}-C_1A_I$ kills
the normal-form span and has rank four. Thus
$\det\sigma_1(L_c)(\xi)=\det(\pi(\xi\odot\,\cdot))$ is the
characteristic quartic up to a nonzero factor.

Apply the linear construction in the proof of
\cite[Theorem~E, Section~5 and Lemma~3.3]{NakamuraMaeda}.
For a literal realization as a linearized first-order
differential map, use $(c,X)\mapsto(E(\partial)c,L_cX)$, where
$E(\xi)$ acts by left quaternionic multiplication on each of
the ten four-vector columns of $c$. Its derivative at $(c,0)$
is block diagonal and has principal determinant
$|\xi|^{40}\det\sigma_1(L_c)$. The normalized scalar symbol
$q=\det\sigma_1(L_c)/|\xi|^3$ has uniformly nonzero $\xi$-gradient
on the nearby cospheres by Euler's identity and regularity.
The normalization only reparametrizes the null bicharacteristics.
On a small outer ball, $(|x|^2)''\ge c_*>0$ along the flow
gives pseudoconvexity. The velocity lower bound and uniform
acceleration bound give a short no-return interval on a smaller ball.
These properties hold in a fixed low Sobolev neighborhood.
For real $f$, take $T_cf$ to be the real part of the last four
components of the inverse (3.15) applied to $(0,f)$. Since $L_c$
is real, $L_cT_cf=f$ on one fixed ball $U$. With an integer $d\ge3$
strictly above the loss and lower threshold in (3.16), the same
inverse satisfies
\[
 \|T_cf\|_{s-d}\le C_s\big(\|f\|_s+(1+\|c\|_s)\|f\|_d\big),
 \qquad s\ge d,\quad s\in\mathbb Z.
\]

This choice also gives joint continuity into $L^2$.
The phase, transport amplitudes, and defect symbols in
(5.4)--(5.5) depend continuously on the coefficient derivatives
needed for the kernel bounds.
On the fixed localized supports, the time-cutoff error is
nonstationary and the evolution defect has order $-5$.
Their kernels converge in Hilbert--Schmidt norm by dominated
convergence, so the localized Neumann inverse is continuous
in operator norm. The time-integrated parametrix converges on Schwartz
data, since its phase and amplitudes are fixed off the spatial cutoff.
Its uniform $L^2$ bound then gives strong convergence on $L^2$.
The order-zero normalized cofactor is a finite sum of fixed
Fourier multipliers with continuously varying left coefficients,
yielding continuity of $T_cf$ into $L^2$.

The localized coefficients are smooth functions of $J^2(\vartheta u)$,
with $H^{S-2}$ norms bounded by $C_S(1+\|\vartheta u\|_S)$.
For large fixed $A_0\ge d+6$, product estimates give
$\|f\|_{d+4}\lesssim\|H\|_{d+4}$ and
$\|f\|_{S-2}\lesssim_S\|H\|_S+(1+\|\vartheta u\|_S)\|H\|_{d+4}$.
Apply the system estimate at $s=S-2$ for $X$ in $H^{S-d-2}$
and at $s=d+4$ for $\|X\|_4\lesssim\|H\|_{d+4}$.
In reconstruction, $H^3\hookrightarrow L^\infty$ controls
the term with a high norm on $M_u^{-1}$. This gives
\[
 \|V\|_{S-d-3}\le C_S(\|H\|_S+\|\vartheta u\|_S\|H\|_{d+4}),
 \qquad S\ge d+4.
\]
Reconstruction gives joint continuity into $H^{-1}$, since the
matrix coefficients converge in $C^1$ and $X$ converges in $L^2$.
Set $D=d+4$. Interpolation with the uniform $H^{S-D+1}$ bound
gives continuity into $H^{S-D}$ for $S\ge A_0$, with the right-inverse
identity and estimate passing to finite-regularity pairs.
\end{proof}

\section{An amplitude normal form with Gram constraints}

Work now in dimension four. For a real scalar $t$, constant in space,
and maps $a:\R^4\to\R^4$, $b:\R^4\to\R^6$, put
\[
 P_t(a,b)=2\operatorname{sym}\partial a+\gamma(b)+t\gamma(a).
\]
For $t>0$, $F_t=(x+ta,\sqrt t\,b)$ satisfies
\begin{equation}\label{eq:amplitude-metric}
 \gamma(F_t)=I+tP_t(a,b).
\end{equation}
At zero set $M=Db(I+tDa)^{-1}$ and
$\widehat b_{ij}=b_{ij}-tMa_{ij}$. Define
\[
 \mathcal S_{ij,kl}=\widehat b_{ij}^T(I_6+tMM^T)^{-1}\widehat b_{kl}.
\]
For $t=r^2>0$ this is the normal Hessian Gram tensor of
$rF_t(x/r)$: the normal columns are $\binom{-\sqrt t M^T}{I_6}$.
The formula gives positivity and rank six near a rank-six reference
also for signed $t$. Put $\mathcal R_t(H)=t^{-1}R_{I+tH}(0)$
with $\mathcal R_0(H)=R^{\rm lin}(H)(0)$. For fixed $(a,b)$ the
Gauss identity $\mathcal G(\mathcal S)=\mathcal R_t(P_t(a,b))$
holds for $t>0$ and extends analytically in $t$.

\begin{theorem}\label{thm:prescribed-gram}
Let $S=B^*B$ be a rank-six positive semidefinite Gram tensor with
regular characteristic data, and let $H_2$ be a homogeneous quadratic
symmetric tensor with $R^{\rm lin}(H_2)=\mathcal G(S)$.
There are polynomials $b_0=\frac12B(x,x)$, a cubic $a_0$ and
a fixed ball $U\Subset B_1(0)$ with the following property.
For each integer $\ell\ge0$ and each $H^\ell$ tolerance there are
$k$ and $\epsilon>0$ such that for $H\in H^\infty(U;\Sym_4)
=\bigcap_{j\ge0}H^j(U;\Sym_4)$ and rank-six positive
semidefinite $S'$ satisfying
\[
 |t|+\|H-H_2\|_{H^k(U)}+\|S'-S\|<\epsilon,
 \qquad \mathcal G(S')=\mathcal R_t(H),
\]
there are smooth $a,b$ within the tolerance of $a_0,b_0$
with $P_t(a,b)=H$ and $\mathcal S(t,a,b)=S'$. They can be chosen
continuously into $H^\ell(U)$ for the indicated data norm. For
each small $(t,H-H_2)$ there is a compatible nearby $S'$.
\end{theorem}

\begin{proof}
Let $T$ be the tangent space of the rank-six Gram manifold at $S$,
consisting of forms with zero $K\times K$ block for $K=\ker B$.
Put $\mathcal V=T\cap\ker\mathcal G$ and let $\Pi$ be the orthogonal
projection onto $\mathcal V$. The map $(\mathcal G,\Pi)$ has
injective derivative on $T$. Composing with a linear left inverse of
that derivative and applying the inverse function theorem gives
local injectivity on the Gram manifold.

For homogeneous quadratic $h$ with $R^{\rm lin}(h)=0$, let
$\Gamma_{ij,k}=\frac12(\partial_i h_{jk}+\partial_j h_{ik}-\partial_k h_{ij})$.
The curvature condition makes $\partial_l\Gamma_{ij,k}$ symmetric
in $l,i,j$. Thus $a_k=\Gamma_{ij,k}(x)x_ix_j/6$ has Hessian
$\Gamma_{ij,k}$. The residual $h-2\operatorname{sym}\partial a$
has zero Christoffel array and hence zero first derivatives.
It vanishes by homogeneity. Apply this to $H_2-\gamma(b_0)$,
whose linearized curvature vanishes by the Gauss identity, to obtain
$A_{\mathrm{ref}}=(a_0,b_0)$ with $P_0(A_{\mathrm{ref}})=H_2$.

Lift $z\in\mathcal V$ to $z=B^*\dot B_z+\dot B_z^*B$ linearly.
Such a lift exists since $B\mapsto B^*B$ is a submersion onto
the rank-six Gram manifold, as is seen in the blocks $K^\perp,K$.
Set $\dot b_z=\frac12\dot B_z(x,x)$. The derivative of the Gauss
identity at $t=0$ shows $R^{\rm lin}(D\gamma_{b_0}\dot b_z)=0$,
so the same cubic primitive cancels that metric variation.
Choose a compactly supported $\chi=1$ near zero and write $Wz$
for the lifted variation after multiplication by $\chi$.
Apply Proposition~\ref{prop:metric-inverse} at
$u_0=(\chi x,\chi b_0)$ and shrink $U$ to lie where $\chi=1$.
Set $T_u=R_u\mathcal E_U$ with a fixed Sobolev extension $\mathcal E_U$.

For $A=(a,b)=\chi A_{\mathrm{ref}}+v$ use the map $(t,v)\mapsto(t,Q,J)$, where
\[
 Q=P_t(A)-H_2,\qquad J=\Pi(\mathcal S(t,A)-S).
\]
On $U$ one has $Q_v=D\gamma_u$ for $u=(\chi x+ta,b)$.
At $(t,v)=(0,0)$ we have $Q_vW=0$ and $J_vW=I$.
Choose an integer $A_*\ge\max(A_0+1,D+6,\ell)$ and use the
Nash--Moser theorem \cite[Section~IV(b)]{BGY} on
\[
 E_j=\R\oplus H^{A_*+j}(\R^4;\R^{10}),\qquad
 F_j=(\R\oplus\mathcal V)\oplus H^{A_*+j-1}(U;\Sym_4).
\]
Use Fourier smoothing with fixed extension and restriction on $U$,
and the identity on the finite-dimensional factors. Rescale the $E_j$
norms by a fixed constant to place the low unit ball near zero.
Sobolev product bounds and the smooth finite-jet dependence of $J$
give the first and second derivative estimates of the theorem.

Write $W_u=(I-T_uQ_v)W$, so $Q_vW_u=0$ on $U$.
On the small low ball, $Q_vW$ is small in $H^{A_*-1}$ and $J_vW$
is close to $I$. Since $H^{A_*-1-D}\hookrightarrow C^2$,
$J_vW_u$ remains invertible with bounded inverse. For derivative
data $(\tau,h,k)$ put $q=h-Q_t\tau$, $w=k-J_t\tau$. The
derivative right inverse is $(\tau,\dot v)$ with
\[
 \dot v=T_uq+W_u(J_vW_u)^{-1}(w-J_vT_uq).
\]
The low estimate bounds the finite coefficient by $C\|(\tau,h,k)\|_{F_\alpha}$,
where $\alpha=D+1$. At order $A_*+j-1$, the norms of $Q_t,Q_vW,u$
are at most $C_j(1+\|v\|_{H^{A_*+j}})$. The metric tame estimate
therefore gives
\[
 \|(\tau,\dot v)\|_{E_{j-\alpha}}
 \le C_j\big(\|(\tau,h,k)\|_{F_j}
             +\|(t,v)\|_{E_j}\|(\tau,h,k)\|_{F_\alpha}\big),\qquad j\ge\alpha.
\]
The joint continuity of $T_u$ at these indices gives the same
property for the displayed right inverse, including $j=\alpha$.

The local selection for independent data $(t,H-H_2,z)$
is continuous from a neighborhood of zero in $F_{12\alpha}$ to
$E_\alpha$, and the same selection takes all-order data to
$E_\infty$. Thus $k=A_*+12\alpha-1$ suffices.
Take $z=\Pi(S'-S)$ for a prescribed compatible target.
Local injectivity and the Gauss identity give $\mathcal S(t,A)=S'$.
Taking $z=0$ gives a nearby compatible Gram tensor for any
nearby $(t,H)$. The low ball enforces the tolerance.
\end{proof}

\section{Geometric consequences}

\subsection{Exact prescribed second jet}
For the metric and $B$ in Theorem~\ref{thm:main}, set $S=B^*B$ and
$H_2(x)=\tfrac12D^2g(0)[x,x]$. Since $dg(0)=0$, $R^{\rm lin}(H_2)=R_g(0)$.
Set $t=r^2$ and $H=(g(r\,\cdot)-I)/r^2$. Then $H\to H_2$
at every fixed Sobolev order and $\mathcal R_t(H)=R_g(0)$
by pullback and constant rescaling. Use $S'=S$ in
Theorem~\ref{thm:prescribed-gram} with $\ell=5$. The map $F(x)=rF_t(x/r)$ induces $g$
on $rU$, and $x+ta$ is lower Lipschitz for small $t$.
An ambient rigid motion sends $F(0)$ and $DF(0)$ to the standard
point and tangent frame; $dg(0)=0$ then implies that $D^2F(0)$
is normal. Two surjective factorizations of the same Gram tensor
differ by an orthogonal normal rotation, giving precisely
the desired $B$.

Finally use the regular common space in dimension four of
Theorem~\ref{thm:gram-lifts}, with $R=R_g(0)$, in normalized
coordinates for any smooth metric. This proves the existence
and prescription assertions of Theorem~\ref{thm:main}.

\subsection{Near-affine embeddings}
\begin{theorem}\label{thm:near-affine}
There is a fixed ball $U\Subset B_1(0)$ such that for each integer
$\ell\ge0$ there exist $k,\epsilon>0$ and $C_\ell$ with the following
property. Every smooth metric $g$ near $\overline U$ with
$\delta=\|g-I\|_{H^k(U)}<\epsilon$ admits a smooth isometric
embedding $F=(F^\parallel,F^\perp)$ on $U$ satisfying
\[
 \|F^\parallel-\mathrm{id}\|_{H^\ell(U)}\le C_\ell\delta,\qquad
 \|F^\perp\|_{H^\ell(U)}\le C_\ell\sqrt{\delta}.
\]
The superscripts denote fixed coordinate projections.
\end{theorem}
\begin{proof}
Use Theorem~\ref{thm:prescribed-gram} at $S=S_*$, $H_2=0$ and Sobolev
order $\max\{\ell,4\}$, so the low bound also controls $C^1$.
If $\delta>0$ choose $t=c\delta$
with fixed $c$ large enough that $(g-I)/t$ is in the small
data ball, and decrease the upper bound on $\delta$ to make $t$
small. Use the selection from its proof with $z=0$.
Then $F_t=(x+ta,\sqrt t\,b)$ embeds $g$
by \eqref{eq:amplitude-metric} and the lower Lipschitz
bound for $x+ta$. The two estimates follow from the uniform
bounds on $a,b$. For $\delta=0$ use the plane.
\end{proof}
The square-root order is sharp for $\ell\ge1$. For $\lambda>0$ and
$g_\lambda=I+\lambda x_1^2dx_2^2$ apply
$Lh=\partial_{11}h_{22}+\partial_{22}h_{11}-2\partial_{12}h_{12}$
to $g_\lambda-I=2\operatorname{sym}Da+\gamma(a)+\gamma(b)$,
where $F=(x+a,b)$. The linear term vanishes and $L(g_\lambda-I)=2\lambda$.
Testing against a fixed compactly supported smooth function of integral
one gives $\lambda\lesssim\|Da\|_2^2+\|Db\|_2^2$ by integration by parts.
Thus $\|a\|_{H^1}=O(\lambda)$ forces $\|b\|_{H^1}\gtrsim\sqrt\lambda$.

Taking $H=0$, $S'=S_*$ and small $t>0$ at the flat reference gives a
regular smooth flat embedding by the same exact Gram identity.
General regular-reference openness follows by omitting the finite
outputs. See also \cite[Section~2]{LongV1}.

\subsection{Quantitative radius}
Fix $U,k,\epsilon$ from Theorem~\ref{thm:near-affine}
with $\ell=0$ and use the same spatial scaling.

\begin{proposition}\label{prop:embedding-radius}
For each integer $m\ge1$ there is a constant $c_m>0$ such that a
smooth metric $g$ near $\overline{B_1(0)}$ with
\[
 \partial^\beta(g-I)(0)=0\ (|\beta|<m),\qquad
 \|g\|_{C^{\max(k,m)}(B_1)}\le K,\quad K\ge1
\]
admits a smooth isometric embedding on $B_{c_mK^{-1/m}}(0)$,
where $k$ is the fixed index just chosen.
\end{proposition}
\begin{proof}
Taylor's formula for $|\beta|<m$ and
$r^{|\beta|}\le r^m$ otherwise give
$\|g(r\,\cdot)-I\|_{H^k(U)}\le C_{U,k,m}K r^m$ for $0<r\le1$.
Take $r=a_mK^{-1/m}$ with a fixed sufficiently small $a_m>0$,
and use the same rescaling, with $c_m=a_m\operatorname{rad}(U)$.
\end{proof}

\subsection{Optimality of the dimension}

\begin{proof}
For the optimality assertion in Theorem~\ref{thm:main}, consider maps
into $\R^9$. The assignment of the $k$-jet of the induced
metric to a $(k+1)$-jet of the map is polynomial. Its domain has
dimension $9\binom{k+5}{4}$ and the metric jet space has dimension
$10\binom{k+4}{4}$. For $k\ge36$ the domain has smaller dimension.
By Sard's theorem the image has measure zero. Choose a jet outside
the image with positive definite zeroth term. Its Taylor polynomial
is an analytic metric near zero, and no smooth map into $\R^9$ induces
it there. The same exclusion holds in smaller Euclidean spaces by
linear inclusion. The embedding assertion provides the upper bound.
\end{proof}

\appendix

\section{An explicit regular flat Gram seed}
\label{app:regular-seed}

We use the signed evaluation method of Blekherman \cite[Section~6]{Blekherman}.
A factorization over $\F_7$ establishes characteristic regularity.

\begin{theorem}
There is an explicit rational, completely symmetric positive
semidefinite tensor $S$ of rank six on $E=\Sym_4$ whose characteristic
quartic is nonsingular over $\C$. With $W=\operatorname{im}S$, the
linear Gauss map $\mathcal G:\Sym^2W\to\mathcal R_4$ is surjective
and has kernel $\R S$.
\end{theorem}

\begin{proof}
Identify symmetric matrices with quadratic forms by
$q_Q(v)=v^TQv$ and use the Frobenius pairing.
Put
\[
\begin{aligned}
q_1&=x(x-y-w),&q_2&=y(y-z-w),&q_3&=z(z-x-w),\\
q_4&=w^2-2w(x+y+z)+8(xy+yz+zx).
\end{aligned}
\]
At $w=0$ the first three equations force $x=y=z$. At $w=1$,
three nonzero coordinates would give $x-y=y-z=z-x=1$.
A single nonzero coordinate equals $1$, while two nonzero
coordinates give a cyclic permutation of $(2,1,0)$.
Thus the common projective zeros over $\C$ are
\[
\begin{gathered}
(0,0,0,1),\ (0,0,1,1),\ (0,1,0,1),\ (0,2,1,1),\\
(1,0,0,1),\ (1,0,2,1),\ (2,1,0,1),\ (1,1,1,0).
\end{gathered}
\]
Denote these representatives by $v_i$ in the displayed order.
The quadratic evaluation map has image $c^\perp$ for
$c=(6,-3,-3,1,-3,1,1,-2)$. By cyclic symmetry the relation
is checked on $w^2,xw,x^2,xy$, with weighted sums
$6-9+3$, $-3+2+1$, $-3+4+1-2$, and $2-2$, respectively.
For the rank, $w^2,xw,yw,zw$ are independent on the four
points with at most one nonzero spatial coordinate and $w=1$.
The mixed monomials $xy,yz,zx$ vanish on those four points
and give a permutation of $2I_3$ on the other points with $w=1$.
Put
\[
 z^*=(q_4(v_i))_i=(1,-1,-1,11,-1,11,11,24),\qquad
 A_i=c_i/z_i^*,\qquad S=\sum_i A_i v_i^{\otimes4}.
\]
Only $A_8$ is negative. For weights of the form $A_i=c_i/z_i^*$
with $z^*\in c^\perp$, one negative weight $A_j$ and the other
weights positive give a positive semidefinite form
$\sum_i A_i z_i^2$ on $c^\perp$ with null line $\R z^*$.
Indeed, eliminating $z_j$
subtracts $(\sum_{i\ne j}c_i z_i)^2/\sum_{i\ne j}c_i^2/A_i$
from $\sum_{i\ne j}A_i z_i^2$, since $c\cdot z^*=0$.
Cauchy--Schwarz gives positivity and the equality case.
Thus $S$ has rank six and kernel $\operatorname{span}(q_i)_{i=1}^4$.
The four quadrics have no common projective zero, so the
complete-intersection calculation in Theorem~\ref{thm:gram-lifts}
gives the claimed Gauss kernel and surjectivity.

For regularity put $p=\det D(q_1,q_2,q_3,q_4)=16p_W$.
Over $\F_7$ use $a=x+y+z$, $b=x+2y+4z$, $c'=x+4y+2z$.
The inverse is obtained by interchanging $2$ and $4$ in the
change matrix and dividing by $3$. Apply the same matrix
to $(q_1,q_2,q_3)$ so that its determinant cancels the input
Jacobian factor. Substituting the inverse and subtracting
twice the new first quadric from the fourth gives
\[
 bc'-aw,\quad -(a+w)b+2c'^2,\quad
 (2a-w)c'-b^2,\quad w^2-2a^2.
\]
Write $A=a+w$, $D=2a-w$. The last Jacobian row is
$(-4a,0,0,2w)$, so its expansion gives
\begin{align*}
 \bar p
 &=4a(aAD+2b^3-4c'^3-2wbc')
   +2w(wAD+2b^3+8c'^3+4abc')\\
 &=f+b^3h+c'^3k,\qquad
 f=(a^2-3w^2)(a+w)(a+3w),\quad
 h=a-3w,\quad k=-2(a-w).
\end{align*}
The binary roots $-1,-3,\pm\sqrt{3}$ of $f(a,1)$ are distinct.
The values $f(3,1)=4$ and $f(1,1)=5$ exclude the zeros of $h$ and $k$.
At a gradient zero over $\overline{\F}_7$, one has
$b^2h=c'^2k=0$ and $\bar p=f=0$ by Euler's identity.
If $(a,w)\ne(0,0)$, then $h,k$ are nonzero, so $b=c'=0$,
contradicting the squarefreeness of $f$. If $(a,w)=(0,0)$,
the binary gradient $b^3\nabla h+c'^3\nabla k$ vanishes
only when $b=c'=0$, since the two linear forms are independent.

Thus the four cubic derivatives form a homogeneous regular
sequence modulo $7$, with quotient Hilbert series $(1+s+s^2)^4$.
Thus the integral multiplication map $(g_i)_i\mapsto\sum_i g_i\partial_i p$
from four copies of degree six to degree nine is onto modulo $7$,
and hence onto over $\C$. At a common complex zero of the
derivatives every ninth power of a coordinate would vanish,
so the characteristic quartic is nonsingular.
\end{proof}

After normalizing $S$ to trace one, the trace-zero lift in
\eqref{eq:gram-selection} is unique since the Gauss kernel is $\R S$.

\section{A block seed for the open-family method}
\label{app:block-seed}

The three quadrics $xy-zw,xz-yw,xw-yz$ have as joint projective
zeros the four coordinate points $e_i$ and the four points
$v_{\epsilon,\eta}=(1,\epsilon,\eta,\epsilon\eta)$,
$\epsilon,\eta\in\{\pm1\}$. Indeed, if a coordinate vanishes
at most one is nonzero, and otherwise setting $x=1$ gives
$y^2=z^2=w^2=1$ and $w=yz$. Quadratic evaluation has rank seven
by the four square columns and three independent sign characters
in the mixed columns. In this order it has relation
$c=(-4,-4,-4,-4,1,1,1,1)$, and $4x^2-y^2-z^2-w^2$ evaluates
to $z^*=(4,-1,-1,-1,1,1,1,1)$. Thus the same signed construction,
divided by four, gives the positive semidefinite rank-six tensor
\[
 S_{\rm B}=\tfrac14\sum_{\epsilon,\eta}v_{\epsilon,\eta}^{\otimes4}
                 +\sum_{i=2}^4e_i^{\otimes4}-\tfrac14e_1^{\otimes4}.
\]
The four kernel quadrics form a regular sequence by the nonzero
evaluations of the fourth. The Gauss kernel is again one-dimensional.
The kernel Jacobian has rank two at $(0,1,1,0)$, so the characteristic
quartic is singular there. The open-family argument chooses
a regular support nearby.

\section{Analytic Gram-jet realization}
\label{app:analytic-jet}

\begin{theorem}
Let $W$ and $R\mapsto S_R$ be the four-dimensional support and
polynomial selection of Theorem~\ref{thm:gram-lifts} or
Appendix~\ref{app:regular-seed}. Every analytic metric $g$ in
coordinates with $g(0)=I$ and $dg(0)=0$ has an analytic local
isometric embedding into $\R^{10}$ whose normal second fundamental
form has Gram tensor $S_{R_g(0)}$ in those coordinates at zero.
\end{theorem}

\begin{proof}
Factor $S_{R_g(0)}=B^*B$ in the original coordinates.
Choose an analytic orthonormal frame $E_i$ whose last vector at
zero is dual to a noncharacteristic covector. The six vectors
$B(E_i(0)\odot E_j(0))$, $1\le i\le j\le3$, are independent by
the slice argument in Proposition~\ref{prop:metric-inverse}.
Prescribe the mixed connection forms using
$h_{ij}=B(E_i(0)\odot E_j(0))$ in the isometric-embedding system.
The Gauss equation and this independence give an ordinary integral
element by the polar calculation in \cite[Ch.~III, proof of
Theorem~3.11]{BCGGG}. Cartan--K\"ahler
\cite[Ch.~III, Corollary~2.3]{BCGGG} gives an integral manifold
with this tangent plane. Its projection
\cite[Ch.~III, Proposition~3.9]{BCGGG} realizes the prescribed
$B$ and hence the Gram tensor in the original coordinates.
The tangential Hessian vanishes since $dg(0)=0$.
\end{proof}

\section*{Declaration of AI use}
The author used OpenAI models as assistive tools in preparing this manuscript.


\begingroup
\small
\begin{thebibliography}{10}
\setlength{\itemsep}{0pt}
\setlength{\parsep}{0pt}
\bibitem{Janet}
M. Janet,
\emph{Sur la possibilit\'e de plonger un espace riemannien donn\'e
dans un espace euclidien},
Ann. Soc. Polon. Math. \textbf{5} (1926), 38--43.
\href{https://jbc.bj.uj.edu.pl/Content/932771/NDIGCZAS094399_1926.pdf}{Original volume, pp.~38--43}.

\bibitem{Cartan}
\'E. Cartan,
\emph{Sur la possibilit\'e de plonger un espace riemannien donn\'e
dans un espace euclidien},
Ann. Soc. Polon. Math. \textbf{6} (1927), 1--7.
\href{https://mccuan.math.gatech.edu/courses/6455/Cartan1927.pdf}{Original volume, pp.~1--7}.

\bibitem{BGY}
R. L. Bryant, P. A. Griffiths, and D. Yang,
\emph{Characteristics and existence of isometric embeddings},
Duke Math. J. \textbf{50} (1983), no.~4, 893--994.
\href{https://doi.org/10.1215/S0012-7094-83-05040-8}{DOI}.

\bibitem{BCGGG}
R. L. Bryant, S. S. Chern, R. B. Gardner, H. L. Goldschmidt,
and P. A. Griffiths,
\emph{Exterior Differential Systems}, Mathematical Sciences
Research Institute Publications, vol.~18, Springer-Verlag,
New York, 1991.
\href{https://doi.org/10.1007/978-1-4613-9714-4}{DOI}.

\bibitem{NakamuraMaeda}
G. Nakamura and Y. Maeda,
\emph{Local smooth isometric embeddings of low-dimensional Riemannian
manifolds into Euclidean spaces},
Trans. Amer. Math. Soc. \textbf{313} (1989), no.~1, 1--51.
\href{https://doi.org/10.1090/S0002-9947-1989-0992597-8}{DOI}.

\bibitem{GoodmanYang}
J. B. Goodman and D. Yang,
\emph{Local solvability of nonlinear partial differential equations
of real principal type}, unpublished manuscript, revised October 11, 1988.
\href{https://cims.nyu.edu/~yangd/papers/goodman-yang.pdf}{Author manuscript}.

\bibitem{HanKhuri}
Q. Han and M. Khuri,
\emph{The linearized system for isometric embeddings and its
characteristic variety},
Adv. Math. \textbf{230} (2012), no.~1, 263--293.
\href{https://doi.org/10.1016/j.aim.2011.11.011}{DOI};
\href{https://arxiv.org/abs/1111.4722}{arXiv:1111.4722}.

\bibitem{Blekherman}
G. Blekherman,
\emph{Nonnegative polynomials and sums of squares},
J. Amer. Math. Soc. \textbf{25} (2012), no.~3, 617--635.
\href{https://doi.org/10.1090/S0894-0347-2012-00733-4}{DOI}.

\bibitem{LongV1}
Z. Long,
\emph{Smooth local isometric embeddings in dimension four},
earlier version, \href{https://arxiv.org/abs/2611.10598v1}{arXiv:2611.10598v1}
(2026).

\bibitem{YangWang}
W. Yang and Y. Wang,
\emph{Smooth local isometric and conformal embeddings in optimal
dimension}, preprint, hexagon:2610.00378v4 (2026).
\href{https://hexagonmath.org/2610.00378v4}{Version~4}.
\end{thebibliography}
\endgroup
\end{document}